\documentclass[12pt,a4paper]{article}
\usepackage[mathcal]{euscript} 
\usepackage{amsmath}
\usepackage{amsthm} 
\usepackage{amssymb}          
\usepackage[T1]{fontenc}
\usepackage[english]{babel}
\usepackage[dvips]{graphicx}
                 
\newcommand{\calC}{\mathcal{C}} 
\newcommand{\calD}{\mathcal{D}} 
\newcommand{\calG}{\mathcal{G}} 
\newcommand{\calM}{\mathcal{M}} 
\newcommand{\calT}{\mathcal{T}} 
\newcommand{\Nset}{\mathbb{N}} 

\newcommand{\leqp}{{\,\preccurlyeq\,}} 
\newcommand{\eqp}{{\,\cong\,}}
\newcommand{\geqp}{{\,\succcurlyeq\,}}
\newcommand{\gp}{{\,\succ\,}}
\newcommand{\lp}{{\,\prec\,}}
\newcommand{\neqp}{{\,\ncong\,}}
\newcommand{\eqs}{{\,\thicksim\,}}
\newcommand{\eqt}{{\,\thickapprox\,}}

\newcommand{\bloks}[2]{{\langle {#1}\rangle_\thicksim^{#2}}}
\newcommand{\blokt}[2]{{\langle {#1}\rangle_\thickapprox^{#2}}}
 
\newcommand{\ie}{i.e.,\ }
\newtheorem{lem}{Lemma}
\newtheorem{thm}[lem]{Theorem} 
\newtheorem{cor}[lem]{Corollary}

\title{Dyck paths and pattern-avoiding matchings}
\author{V\'\i t Jel\'\i nek\\
\small Institute for Theoretical Computer Science (ITI)\thanks{ITI is supported by project 
1M0021620808 of the Ministry of Education of the Czech Republic.},\\ 
\small Charles University\\
\small Malostransk\'e n\'am\v est\'\i~25\\
\small Prague, Czech Republic.\\
\small E-mail:\texttt{ jelinek@kam.mff.cuni.cz}
} 
\date{}

\begin{document}

\maketitle

\begin{abstract}
How many matchings on the vertex set $V=\{1,2,\dotsc,2n\}$  avoid a given 
configuration of three edges? Chen, Deng and Du have shown that the number of 
matchings that avoid three nesting edges is equal to the number of matchings 
avoiding three pairwise crossing edges. In this paper, we consider other 
forbidden configurations of size three. We present a bijection between 
matchings avoiding three crossing edges and matchings avoiding an edge nested 
below two crossing edges. This bijection uses non-crossing pairs of Dyck paths 
of length $2n$ as an intermediate step. 

Apart from that, we give a bijection that maps matchings avoiding two nested 
edges crossed by a third edge onto the matchings avoiding all configurations 
from an infinite family $\mathcal{M}$, which contains the configuration 
consisting of three crossing edges. We use this bijection to show that for 
matchings of size $n>3$, it is easier to avoid three crossing edges than to 
avoid two nested edges crossed by a third edge.

 In this updated version of this paper, we add new references 
to papers that have obtained analogous results in a different context.
\end{abstract}

\section{Introduction and Basic Definitions}
                   
This is an updated preprint of a paper whose journal version has 
already appeared in print \cite{jour}. The main reason for the 
update was to include references to the papers \cite{bwx,s,sw}, which have 
independently obtained equivalent results using different methods.

The enumeration of pattern-avoiding permutations has received a considerable 
amount of attention lately (see \cite{kima} for a survey). We say that a 
permutation $\pi$ of order $n$ \emph{contains} a permutation $\sigma$ of 
order $k$, if there is a sequence $1\le i_1<i_2<\dotso<i_k\le n$ such that 
for every $s,t\in[n]$  $\pi(i_s)<\pi(i_t)$ if and only if 
$\sigma(s)<\sigma(t)$. One of the central notions in the study of 
pattern-avoiding permutations is the \emph{Wilf equivalence}: we say that a 
permutation $\sigma_1$ is Wilf-equivalent to a permutation $\sigma_2$ if, for 
every $n\in\Nset$, the number of permutations of order $n$ that avoid 
$\sigma_1$ is equal to the number of permutations of order $n$ that avoid 
$\sigma_2$. In this paper, we consider pattern avoidance in matchings. This 
is a more general concept than pattern avoidance in permutations, since every 
permutation can be represented by a matching. 

A \emph{matching} of size $m$ is a graph on the vertex set 
$[2m]=\{1,2,\dotsc,2m\}$ whose every vertex has degree one. We say that a matching 
$M=(V,E)$ \emph{contains} a matching $M'=(V',E')$ if there is a monotone 
edge-preserving injection from $V'$ to $V$; in other words, $M$ contains $M'$ 
if there is a function $f\colon V'\to V$ such that $u<v$ implies $f(u)<f(v)$ and $\{u,v\}\in 
E'$ implies $\{f(u),f(v)\}\in E$.

Let $M$ be a matching of size $m$, and let $e=\{i,j\}$ be an arbitrary edge of~$M$. 
If $i<j$, we say that $i$ is an \emph{l-vertex} and $j$ is an 
\emph{r-vertex} of $M$. Obviously, $M$ has $m$ l-vertices and $m$ r-vertices. 
Let $e_1=\{i_1,j_1\}$ and $e_2=\{i_2,j_2\}$ be two edges of $M$, with 
$i_1<j_1$ and $i_2<j_2$, and assume that $i_1<i_2$. We say that the two edges 
$e_1$ and $e_2$ \emph{cross} each other if $i_1<i_2<j_1<j_2$, and we say that 
$e_2$ is \emph{nested} below $e_1$ if $i_1<i_2<j_2<j_1$.                                                      

\begin{figure}
\begin{tabular}{ccc}
\includegraphics[scale=.65]{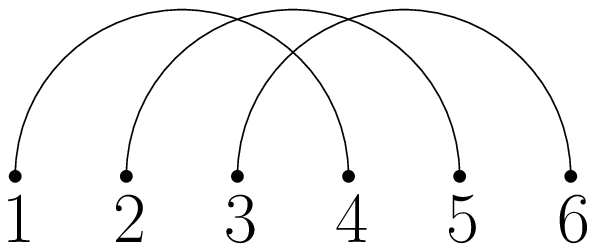}
&
\includegraphics[scale=.65]{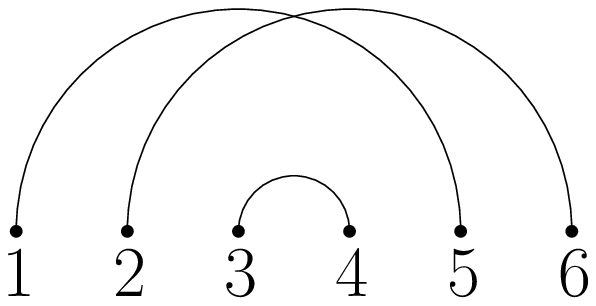}
&
\includegraphics[scale=.65]{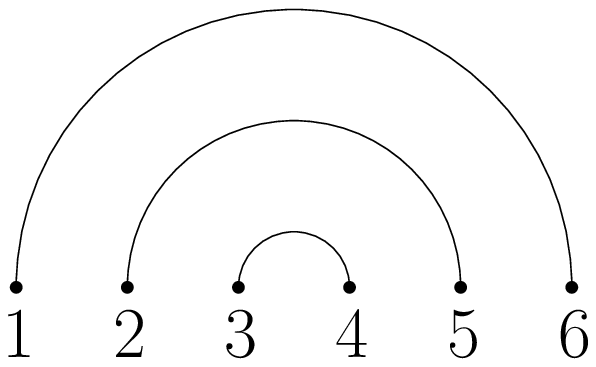}\\
$M_{123}$&$M_{231}$&$M_{321}$\\[10pt]
\includegraphics[scale=.65]{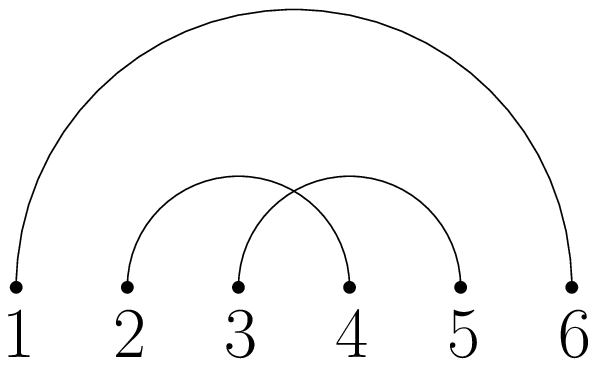}
&
\includegraphics[scale=.65]{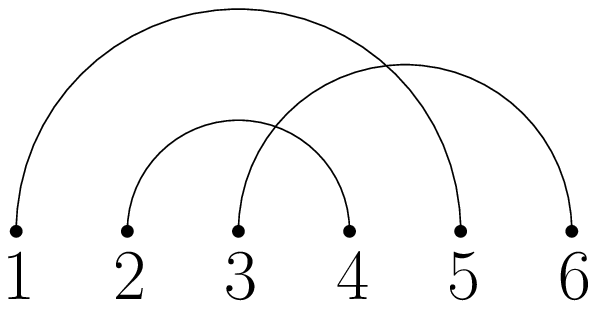}
&
\includegraphics[scale=.65]{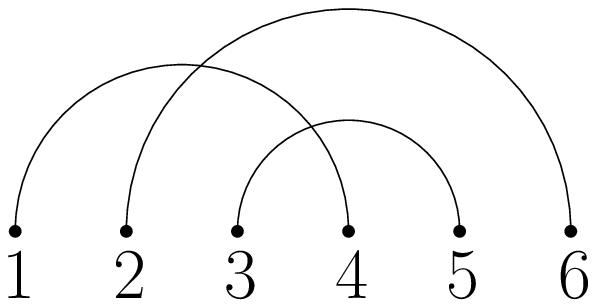}\\
$M_{312}$&$M_{213}$&$M_{132}$\\
\end{tabular}
\caption[Permutational matchings with three edges]{The 
six permutational matchings with three edges.}\label{fig-match}
\end{figure}

We say that a matching $M$ on the vertex set $[2m]$ is \emph{permutational}, 
if for every l-vertex $i$ and every r-vertex $j$ we have $i\le 
m<j$. There is a natural one-to-one correspondence between permutations of 
order $m$ and permutational matchings of size $m$: if $\pi$ is a 
permutation of $m$ elements, we let $M_\pi$ denote the permutational 
matching on the vertex set $[2m]$ whose edge set is the set 
$\bigr\{\{i,m+\pi(i)\},\ i\in[m]\bigl\}$. In this paper, we will often 
represent a permutation $\pi$ on $m$ elements by the ordered sequence 
$\pi(1)\pi(2)\dotsm\pi(m)$. Thus, for instance, $M_{132}$ refers to the 
matching on the vertex set $[6]$, with edge set 
$\bigr\{\{1,4\},\{2,6\},\{3,5\}\bigl\}$. Figure~\ref{fig-match} depicts most 
of the matchings relevant for this paper. Note that a permutational 
matching $M_\pi$ contains the permutational matching $M_\sigma$ if and 
only if $\pi$ contains~$\sigma$.
                                   
Let $n=2m$ be an even number. A \emph{Dyck path} of length $n$ is a piecewise 
linear nonnegative walk in the plane, which starts at the 
point $(0,0)$, ends at the point $(n,0)$, and consists of $n$ linear segments 
(``steps''), of which there are two kinds: an \emph{up-step} connects $(x,y)$ 
with $(x+1,y+1)$, whereas a \emph{down-step} connects $(x,y)$ with 
$(x+1,y-1)$. The nonnegativity of the path implies that among the first $k$ 
steps of the path there are at least $k/2$ up-steps. Let $\calD_m$ denote the 
set of all Dyck paths of length $2m$. It is well known that 
$|\calD_{m}|=c_m$, where $c_m=\frac{1}{m+1}\binom{2m}{m}$ is the $m$-th 
Catalan number. 

Every Dyck path $D\in\calD_m$ can be represented by a \emph{Dyck word} 
(denoted by $w(D)$), which is a binary word $w\in\{0,1\}^{2m}$ such that 
$w_i=0$ if the $i$-th step of $D$ is an up-step, and $w_i=1$ if the 
$i$-th step of $D$ is a down-step. It can be easily seen that a word 
$w\in\{0,1\}^n$ is a Dyck word of some Dyck path if and only if the 
following conditions are satisfied:
\begin{itemize}
\item
The length $n=|w|$ is even.
\item 
The word $w$ has exactly $n/2$ terms equal to 1.
\item
Every prefix $w'$ of $w$ has at most $|w'|/2$ terms equal to 1.
\end{itemize}
We will use the term Dyck word to refer to any binary word satisfying 
these conditions. The set of all Dyck words of length $2m$ will be 
denoted by~$\calD'_m$.

Let $\calG(m)$ denote the set of all matchings on the vertex set 
$[2m]$. For a matching $M\in\calG(m)$, we define the \emph{base} of $M$ 
(denoted by $b(M)$) to be the binary word $w\in\{0,1\}^{2m}$ such that 
$w_i=0$ if $i$ is an l-vertex of $M$, and $w_i=1$ if $i$ is an r-vertex of 
$M$. The base $b(M)$ is clearly a Dyck word; conversely, every Dyck word is a 
base of some matching. If $w_i=0$ (or $w_i=1$) we say that $i$ is an 
\emph{l-vertex} (or an \emph{r-vertex}, respectively) \emph{with respect to 
the base $w$}. Let $m\in\Nset$, let $\calM$ be an arbitrary set of 
matchings, and let $w\in\calD'_m$; we define the sets $\calG(m,\calM)$ and 
$\calG(m,w,\calM)$ as follows:
\begin{align*}
\calG(m,\calM)&=\{M\in\calG(m);\ M\text{ avoids all the elements of }\calM\}\\
\calG(m,w,\calM)&=\{M\in\calG(m,\calM);\ b(M)=w\}
\end{align*}

Let $g(m), g(m,\calM)$ and $g(m,w,\calM)$ denote the cardinalities of 
the sets $\calG(m)$, $\calG(m,\calM)$ and $\calG(m,w,\calM)$, 
respectively. The sets $\calG(m,w,\calM)$ form a partition of 
$\calG(m,\calM)$. In other words, we have
\[
\calG(m,\calM)=\bigcup_w \calG(m,w,\calM)\quad
\text{ and }\quad g(m,\calM)=\sum_w g(m,w,\calM),
\]
where the union and the sum range over all Dyck words $w\in\calD'_m$.

If no confusion can arise, we will write $\calG(m,M)$ 
instead of $\calG(m,\{M\})$ and $\calG(m,w,M)$ instead of $\calG(m,w,\{M\})$.

 There is an alternative way to encode ordered matchings, which 
uses transversals of Ferrers shapes. A \emph{Ferrers shape} is a 
left-justified array of cells, where the number of cells in a given row does 
not exceed the number of cells in the row directly above it. A 
\emph{transversal} of a Ferrers shape is a subset of its cells which 
intersects every row and every column exactly once. A matching $M$ on the 
vertex set $[2m]$ can be represented by a transversal of a Ferrers shape as 
follows: first, consider a Ferrers shape with $m$ rows and $m$ columns, in 
which the $i$-th row (counted from the bottom) has $k$ cells if and only if 
the $i$-th r-vertex of $M$ has $k$ l-vertices to the left of it. Next, we 
define a transversal of this shape: the transversal contains the cell in row 
$i$ and column $j$ if and only if the $i$-th r-vertex is connected to the 
$j$-th l-vertex by an edge of $M$. This correspondence establishes a 
bijection between transversals and matchings, where matchings of a given base 
correspond to the transversals of a given shape. In the context of 
pattern-avoiding transversals, the equivalence relation $\eqp$ is known as 
shape-Wilf equivalence. Note that the permutational matchings correspond to 
transversals of square shapes, i.e., to permutation matrices.

The correspondence between matchings and transversals, as well as the 
(more general) correspondence between ordered graphs and nonnegative 
fillings of Ferrers shapes, has been pointed out by Krattenthaler~ \cite{kra} and 
de Mier~\cite{adm}. Using this correspondence, it becomes clear that some previous 
results on pattern-avoiding graphs \cite{chdd,stan,du,jour} are equivalent to 
results on pattern-avoiding fillings \cite{bwx,kra,adm,s,sw}.
                                                       
The aim of this paper is to study the relative cardinalities of the sets 
$\calG(m,F)$, with $F$ being a permutational matching with three edges. 
For this purpose, we introduce the following notation:

Let $\leqp$ be the quasiorder relation defined as follows: for two sets 
$\calM$ and $\calM'$ of matchings, we write $\calM\leqp\calM'$, 
if for each $m\in\Nset$ and each $w\in\calD'_m$ we have $g(m,w,\calM)\le 
g(m,w,\calM')$. Similarly, we write $\calM\eqp \calM'$ if $\calM\leqp 
\calM'$ and $\calM\geqp \calM'$, and we write $\calM\lp \calM'$ if 
$\calM\leqp \calM'$ and $\calM\neqp \calM'$. As above, we omit the curly 
braces when the arguments of these relations are singleton sets.
                                                           
Note that two permutations $\pi$ and $\sigma$ are Wilf-equivalent if and only 
if for every $m\in\Nset$ the equality $g(m,0^m1^m,M_\sigma)=g(m,0^m1^m,M_\pi)$ holds, 
where $0^m1^m$ is the Dyck word consisting 
of $m$ consecutive $0$-terms followed by $m$ consecutive $1$-terms. Thus, if 
$M_\pi\eqp M_\sigma$, then $\pi$ and $\sigma$ are Wilf-equivalent; however, 
the converse does not hold in general: it is well known that all the 
permutations of order three are Wilf-equivalent, whereas the results of this 
paper imply that the permutational matchings of size three fall into three $\eqp$-classes.

Combining the known results on Ferrers transversals and the 
known results on matchings, the full characterization of the $\eqp$ 
and $\lp$ relations for patterns of size three has been obtained:

\[ 
M_{213}\eqp M_{132} \lp M_{123}\eqp M_{321} \eqp M_{231}\lp M_{312}.
\]         

The equivalence $M_{213}\eqp M_{132}$ follows from the results of Stankova 
and West~\cite{sw} obtained in the context of fillings of Ferrers shapes. 
In~\cite{jour}, the same equivalence is proved in the context of 
pattern-avoiding matchings as a corollary to the result presented in this 
paper as Theorem~\ref{thm-bij}. The equivalence of $M_{123}\eqp M_{321}$ 
follows from the more general result $M_{12\dotsb k}\eqp M_{k(k-1)\dotsb 1}$, 
which was proved by Chen et al.~\cite{chdd} in the context of matchings, and 
by Backelin et al.~\cite{bwx} in the context of transversals. Several 
generalizations of this result are obtained by Krattenthaler~\cite{kra} (see 
also~\cite{adm}). The equivalence $M_{321} \eqp M_{231}$ follows from the 
general results of Backelin et al.~\cite{bwx}. In this paper, a different 
bijective argument is given (see Section~\ref{sec-231}). The relation 
$M_{132} \lp M_{123}$ is proved in Corollary~\ref{cor-123-132} of this paper, 
and a different argument is given in a forthcoming paper by 
Stankova~\cite{s}, which also contains the proof of $M_{231}\lp M_{312}$, 
completing the classification. 
   
This paper is organized as follows: in Section~\ref{sec-132}, we prove 
Theorem~\ref{thm-bij}, which simultaneously implies $M_{213}\eqp M_{132}$ and 
$M_{132} \lp M_{123}$. In Section~\ref{sec-231}, we present a bijective 
argument that implies $M_{321} \eqp M_{231}$.  

\section{The forbidden matchings $M_{132}$ and $M_{213}$}\label{sec-132}

Since the matching $M_{132}$ is the mirror image of the matching 
$M_{213}$, it is obvious that $g(m,M_{132})$ is equal to $g(m,M_{213})$ 
for each $m\in \Nset$. However, there seems to be no straightforward 
argument demonstrating the stronger fact that $M_{132}\eqp M_{213}$.

For $k\ge3$, we define $C_k\in\calG(k)$ to be the matching with 
edge set $E(C_k)=\bigl\{\{2i-1,2i+2\};\ 1\le i<k\bigr\}\cup\bigl\{\{2,2k-1\}\bigr\}$. 
Let $\calC=\{C_k;\ k\ge 3\}$. 

The goal of this section is to prove the following result:
\begin{thm}\label{thm-bij}
$\calC\eqp M_{132}$.
\end{thm} 

Since all the 
elements of $\calC$ are symmetric upon mirror reflection, this also proves 
that $\calC\eqp M_{213}$ and $M_{132}\eqp M_{213}$, see 
Corollaries~\ref{cor-calc} and~\ref{cor-sym} at the end of this section. 

Throughout this section, we consider $m\in\Nset$ and $w\in\calD'_m$ to be 
arbitrary but fixed, and we let $n=2m$. For the sake of brevity, we 
write $\calG^M$ instead of $\calG(m,w,M_{132})$ and $\calG^C$ 
instead of $\calG(m,w,\calC)$. For a matching $G\in\calG(m)$ and an 
arbitrary integer $k\in[n]$, let $G[k]$ denote the subgraph of $G$ 
induced by the vertices in $[k]$. There are three types of vertices in $G[k]$: 
\begin{itemize}
\item
The r-vertices of $G$ belonging to $[k]$. Clearly, all these vertices 
have degree one in $G[k]$.
\item
The l-vertices of $G$ connected to some r-vertex belonging to $[k]$. 
These have degree one in $G[k]$ as well.
\item
The l-vertices of $G$ belonging to $[k]$ but not connected to an 
r-vertex belonging to $[k]$. These are the isolated vertices of $G[k]$, 
and we will refer to them as the \emph{stubs} of $G[k]$.
\end{itemize}

Let $G$ be an arbitrary graph from $\calG^{M}$. The sequence 
\[
G[1],\ G[2],\ G[3], \dotsc,\ G[n-1],\ G[n]=G 
\]
will be called \emph{the construction} of $G$. It is convenient to view 
the construction of $G$ as a sequence of steps of an algorithm that 
produces the matching $G$ by adding one vertex in every step. Two graphs 
$G, G'$ from $\calG^{M}$ may share an initial part of their 
construction; however, if $G[k]\neq G'[k]$ for some $k$, then obviously 
$G[j]\neq G'[j]$ for every $j\ge k$. It is natural to represent the set 
of all the constructions of graphs from $\calG^{M}$ by \emph{the 
generating tree} of $\calG^{M}$ (denoted by $\calT^{M}$), defined by 
the following properties:
\begin{itemize}
\item
The generating tree is a rooted tree with $n$ levels, where the root 
is the only node of level one, and all the leaves appear on level $n$.
\item
The nodes of the tree are exactly the elements of the following set: 
\[
\{G';\ \exists k\in[n], \exists G\in\calG^{M}:\  G'=G[k]\}
\]
\item
The children of a node $G'$ are exactly the elements of the following set:
\[
\{G''; \exists k\in[n-1], \exists G\in\calG^{M}:\ G'=G[k], G''=G[k+1]\}
\]
\end{itemize}
It follows that the level of every node $G'$ of the tree $\calT^M$ is 
equal to the number of vertices of $G'$. Also, the leaves of $\calT^M$ 
are exactly the elements of $\calG^{M}$, and the nodes of the path from 
the root to a leaf $G$ form the construction of $G$.

The generating tree of $\calG^{C}$, denoted by $\calT^{C}$, is defined 
in complete analogy with the tree $\calT^M$. Our goal is 
to prove that the two trees are isomorphic, hence they have 
the same number of leaves, \ie $|\calG^{M}|=|\calG^{C}|$.

We say that a graph $G'$ on the vertex set $[k]$ is 
\emph{consistent} with $w$, if $G'=G[k]$ for some matching $G\in 
\calG(m)$ with base $w$.

\begin{lem}\label{lem-t1t2}\ 

\noindent
1. A graph $G'$ is a node of\/ $\calT^{M}$ if and only if $G'$ satisfies these three conditions:
\begin{itemize}
\item[(a)]
$G'$ is consistent with $w$. 
\item[(b)]
$G'$ avoids $M_{132}$.
\item[(c)]
$G'$ does not contain a sequence of five vertices $x_1<x_2<\dotsb<x_5$ 
such that $x_2$ is a stub, while $\{x_1,x_4\}$ and $\{x_3,x_5\}$ are 
edges of $G'$.
\end{itemize}

\noindent 2. A graph $H'$ is a node of\/ $\calT^{C}$ if and only if $H'$ satisfies these three conditions:
\begin{itemize}
\item[(a)]
$H'$ is consistent with $w$. 
\item[(b)]
$H'$ avoids $\calC$.
\item[(c)]                                     
For every $p\ge 3$, $H'$ does not contain an induced subgraph isomorphic to 
$C_p[2p-1]$ (by an order preserving isomorphism). In other words, for every 
$p\ge 3$, $H'$ does not contain a sequence of\/ $2p-1$ vertices 
$x_1<x_2<\dotsb<x_{2p-1}$, where $2p-3$ is a stub, and the remaining $2p-2$ 
vertices induce the edges $\bigl\{\{x_{2i-1},x_{2i+2}\}, 1\le i\le 
p-2\bigr\}\cup\bigl\{\{x_2,x_{2p-1}\}\bigr\}$. 
\end{itemize}
\end{lem}
\begin{proof}
We first prove the first part of the lemma. Let $G'$ be a node of 
$\calT^{M}$. Clearly, $G'$ satisfies conditions \textit{a} and \textit{b} of 
the first part of the lemma. Assume that $G'$ fails to satisfy condition 
\textit{c}. Choose $G\in\calG^{M}$ such that $G'=G[k]$ for some $k\in[n]$. 
Let $x_6$ denote the r-vertex of $G$ connected to $x_2$. Then $x_6>k$, 
because $x_2$ was a stub of $G'=G[k]$, which implies that $x_6>x_5$ and the 
six vertices $x_1<\dotsb<x_6$ induce a subgraph isomorphic to $M_{132}$, 
which is forbidden. This shows that the conditions \textit{a}, \textit{b} and 
\textit{c} are necessary.

To prove the converse, assume that $G'$ satisfies the three conditions, 
and let $V(G')=[k]$. We will extend $G'$ into a graph $G$ with base $w$, 
by adding the vertices $k+1, k+2,\dotsc, n$ one by one, and each time 
that we add a new r-vertex $i$, we connect $i$ with the smallest stub of 
the graph constructed in the previous steps. We claim that this 
algorithm yields a graph $G\in\calG^{M}$. For contradiction, assume that 
this is not the case, and that there are six vertices 
$x_1<x_2<\dotsb<x_6$ inducing a copy of $M_{132}$. By 
condition \textit{b}, we know that these six vertices are not all 
contained in $G'$, which means that $x_6>k$. Also, by condition 
\textit{c}, we know that $x_5>k$. In the step of the above construction 
when we added the r-vertex $x_5$, both $x_2$ and $x_3$ were stubs. Since 
$x_5$ should have been connected to the smallest available stub, it 
could not have been connected to $x_3$, which contradicts the assumption 
that $x_1,\dotsc,x_6$ induce a copy of $M_{132}$. Thus 
$G\in\calG^{M}$, as claimed.

The proof of the second part of the lemma follows along the same lines. To 
see that the conditions \textit{a}, \textit{b} and \textit{c} of the second 
part are sufficient, note that every graph satisfying these conditions can be 
extended into a graph $H\in\calG^{C}$ by adding new vertices one by one, and 
connecting every new r-vertex to the biggest stub available when the r-vertex 
is added. We omit the details.
\end{proof}

Let $G'$ be a node of $\calT^{M}$. We define a binary relation $\eqs$ on 
the set of stubs of $G'$ by the following rule: $u\eqs v$ if and only if 
either $u=v$ or there is an edge $\{x,y\}\in E(G')$ such that $x<u<y$ 
and $x<v<y$.                                                      

\begin{figure}
\hfil
\includegraphics[scale=.9]{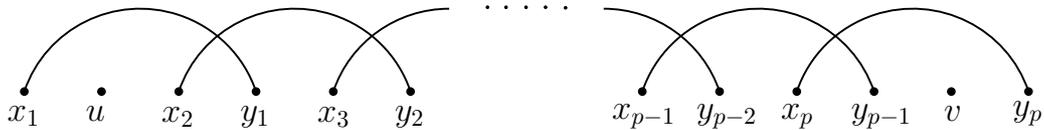}
\caption[The relation $\eqt$.]{The relation $\eqt$. Here $u\eqt v$, 
assuming $u$ and $v$ are stubs.}\label{fig-eqst}
\end{figure} 

Let $H'$ be a node of $\calT^{C}$. We define a binary relation $\eqt$ on the 
set of stubs of $H'$ by the following rule: $u\eqt v$ if and only if either 
$u=v$ or $H'$ contains a sequence of edges $e_1,e_2,\dotsc,e_p$, where $p\ge 
1$, $e_i=\{x_i,y_i\}$, the edge $e_i$ crosses the edge $e_{i+1}$ for each 
$i<p$, and at the same time $x_1<u<y_1$ and $x_p<v<y_p$ (see 
Fig.~\ref{fig-eqst}; note that we may assume, without loss of generality, 
that the edge $e_i$ does not cross any other edge of the sequence except for 
$e_{i-1}$ and $e_{i+1}$, and that no two edges of the sequence are nested: 
indeed, a minimal sequence $(e_i)_{i=1}^p$ witnessing $u\eqt v$ clearly 
has these properties). We remark that the relation $\eqt$ has an intuitive
geometric interpretation: assume that the vertices of $H'$ are represented by 
points on a horizontal line, ordered left-to-right according to the natural 
order, and assume that every edge of $H'$ is represented by a half-circle 
connecting the corresponding endpoints. Then $u\eqt v$ if and only if every vertical 
line separating $u$ from $v$ intersects at least one edge of $H'$.
 
Using condition \textit{c} of the first part of Lemma~\ref{lem-t1t2}, it 
can be easily verified that for every node $G'$ of the tree $\calT^{M}$, 
the relation $\eqs$ is an equivalence relation on the set of stubs of 
$G'$. Let $\bloks{x}{G'}$ denote the block of $\eqs$ containing the stub 
$x$. Clearly, the blocks of $\eqs$ are contiguous with respect to the 
ordering $<$ of the stubs of $G'$; \ie if $x<y<z$ are three 
stubs of $G'$, then $x\eqs z$ implies $x\eqs y\eqs z$.

Similarly, $\eqt$ is an equivalence relation on the set of stubs of a 
node $H'$ of $\calT^{C}$ (notice that, contrary to the case of $\eqs$, 
the fact that $\eqt$ is an equivalence relation does not rely on the 
particular properties of the nodes of $\calT^{C}$ described in 
Lemma~\ref{lem-t1t2}). The block of $\eqt$ containing $x$ will be 
denoted by $\blokt{x}{H'}$. These blocks are contiguous with respect to the 
ordering $<$ as well. 

\begin{lem}\label{lem-child}\ 

\noindent
1. Let $G'$ be a node of level $k<n$ in the tree $\calT^{M}$. The 
following holds:
\begin{itemize}
\item[(a)]
Let $G''$ be an arbitrary child of\/ $G'$ in the tree $\calT^{M}$. This 
implies that $V(G'')=[k+1]$. If the vertex $k+1$ is an l-vertex with 
respect to $w$, then $G''$ is the only child of\/ $G'$, and $k+1$ is a 
stub in $G''$. In this case, $\bloks{x}{G'}=\bloks{x}{G''}$ for every 
stub $x$ of $G'$, and $\bloks{k+1}{G''}=\{k+1\}$. On the other hand, if 
$k+1$ is an r-vertex, then in the graph $G''$ the vertex $k+1$ is 
connected to a vertex $x$ satisfying $x=\min\bloks{x}{G'}$. In this 
case, we have $\bloks{y}{G'}=\bloks{y}{G''}$ whenever $y<x$, and all the 
stubs $z>x$ of $G''$ form a single $\eqs\!$-block in $G''$.
\item[(b)]
If $k+1$ is an r-vertex, then for every stub $x$ satisfying 
$x=\min\bloks{x}{G'}$, $G'$ has a child $G''$ which contains the edge 
$\{x,k+1\}$. This implies, together with part a, that if $k+1$ is an 
r-vertex, then the number of children of $G'$ in $\calT^{M}$ is equal to 
the number of its $\eqs\!$-blocks.
\end{itemize}

\noindent
2. Let $H'$ be a node of level $k<n$ in the tree $\calT^{C}$. The 
following holds:
\begin{itemize}
\item[(a)]
Let $H''$ be an arbitrary child of $H'$ in the tree $\calT^{C}$. This 
implies that $V(H'')=[k+1]$. If the vertex $k+1$ is an l-vertex with 
respect to $w$, then $H''$ is the only child of $H'$, and $k+1$ is a 
stub in $H''$. In this case, $\blokt{x}{H'}=\blokt{x}{H''}$ for every 
stub $x$ of $H'$, and $\blokt{k+1}{H''}=\{k+1\}$. On the other hand, if 
$k+1$ is an r-vertex, then in the graph $H''$ the vertex $k+1$ is 
connected to a vertex $x$ satisfying $x=\max\blokt{x}{H'}$. In this 
case, we have $\blokt{y}{H'}=\blokt{y}{H''}$ whenever $y<x$ and 
$y\not\in\blokt{x}{H'}$, and all the other stubs of $H''$ form a 
single $\eqt\!$-block in $H''$.
\item[(b)]
If $k+1$ is an r-vertex, then for every stub $x$ satisfying 
$x=\max\blokt{x}{H'}$, $H'$ has a child $H''$ which contains the edge 
$\{x,k+1\}$. This implies, together with part a, that if $k+1$ is an 
r-vertex, then the number of children of $H'$ in $\calT^{C}$ is equal to 
the number of its $\eqt\!$-blocks.
\end{itemize}
\end{lem}
\begin{proof}
We first prove part $1a$. The case when $k+1$ is an l-vertex follows 
directly from the definition of $\eqs$, so let us assume that $k+1$ is 
an r-vertex, and let $x$ be the vertex connected to $k+1$ in $G''$. 
Assume, for contradiction, that $x\not=\min\bloks{x}{G'}$, and choose 
$y\in\bloks{x}{G'}$ such that $y<x$. Since $y\eqs x$, $G'$ must contain 
an edge $e=\{u,v\}$, with $u<y<x<v$. Then the five vertices $u, y, x, v, 
k+1$ form in $G''$ a configuration that was forbidden by 
Lemma~\ref{lem-t1t2}, part $1c$. This shows that $x=\min\bloks{x}{G'}$. 
The edge $\{x,k+1\}$ guarantees that all the stubs larger than $x$ are 
$\eqs\!$-equivalent in $G''$, whereas the equivalence classes of the 
stubs smaller than $x$ are unaffected by this edge. This concludes the 
proof of part $1a$.

To prove part $1b$, it is sufficient to show that after choosing a 
vertex $x$ such that $x=\min\bloks{x}{G'}$ and adding the edge 
$\{x,k+1\}$ to $G'$, the resulting graph $G''$ satisfies the three 
conditions of the first part of Lemma~\ref{lem-t1t2}. Condition $1a$ of 
Lemma~\ref{lem-t1t2} is satisfied automatically. If $G''$ fails to 
satisfy condition $1b$, then $G'$ fails to satisfy one of the 
conditions $1b$ and $1c$ of that Lemma, which is impossible. Similarly, if 
$G''$ fails to satisfy condition $1c$, then either $G'$ fails to 
satisfy this condition as well, or $G'$ contains a stub $y$ with 
$y<x$ and $y\eqs x$, contradicting our choice of $x$. 

The proof of the second part of this lemma follows along the same lines 
as the proof of the first part, and we omit it.
\end{proof}

We are now ready to state and prove the main theorem of this section.

\begin{thm}\label{thm-m132} The trees $\calT^{M}$ and $\calT^{C}$ are 
isomorphic. \end{thm} 
\begin{proof}
Our aim is to construct a mapping $\phi$ with the following properties:
\begin{itemize}
\item
The mapping $\phi$ maps the nodes of $\calT^{M}$ to the nodes of 
$\calT^{C}$, preserving their level.
\item
If $G'$ is a child of $G$ in $\calT^M$, then $\phi(G')$ is a child 
$\phi(G)$ in $\calT^C$. Furthermore, if $G_1$ and $G_2$ are two distinct 
children of a node $G$ in $\calT^{M}$, then $\phi(G_1)$ and $\phi(G_2)$ 
are two distinct children of $\phi(G)$ in $\calT^{C}$.
\item
Let $G$ be an arbitrary node of $\calT^{M}$, and let $H=\phi(G)$. Let 
$\bloks{x_1}{G}$, $\bloks{x_2}{G},\, \dotsc,$ $\bloks{x_s}{G}$ be the 
sequence of all the distinct blocks of $\eqs$ in $G$, uniquely 
determined by the condition $x_1<x_2<\dotsb<x_s$. Similarly, let 
$\blokt{y_1}{H}$, $\blokt{y_2}{H},\, \dotsc,$ $\blokt{y_t}{H}$ be the 
sequence of all the distinct blocks of $\eqt$ in $H$, uniquely 
determined by the condition $y_1<y_2<\dotsb<y_t$. Then $s=t$ and 
$|\bloks{x_i}{G}|=|\blokt{y_i}{H}|$ for each $i\in [s]$.
\end{itemize}
These conditions guarantee that $\phi$ is an isomorphism, because, 
thanks to Lemma~\ref{lem-child}, we know that the number of children of 
each node of $\calT^{M}$ (or~$\calT^{C}$) at level $k$ is either equal 
to one if $k+1$ is an l-vertex or equal to the number of blocks of its 
$\eqs$ relation (or $\eqt$ relation, respectively) if $k+1$ is an 
r-vertex.

The mapping $\phi$ is defined recursively for nodes of increasing level. 
The root of $\calT^{M}$ is mapped to the root of $\calT^{C}$. Assume 
that the mapping $\phi$ has been determined for all the nodes of 
$\calT^{M}$ of level at most $k$, for some $k\in[n-1]$, and that it does 
not violate the properties stated above. Let $G$ be a node of level $k$, 
let $H=\phi(G)$. If $k+1$ is an l-vertex, then $G$ has a unique child 
$G'$ and $H$ has a unique child $H'$. In this case, define 
$\phi(G')=H'$. Let us now assume that $k+1$ is an r-vertex. Let 
$\bloks{x_1}{G}$, $\bloks{x_2}{G},\, \dotsc,$ $\bloks{x_s}{G}$ be the 
sequence of all the distinct blocks of $\eqs$ on $G$, with 
$x_1<x_2<\dotsb<x_s$. We may assume, without loss of generality, that 
$x_i=\min\bloks{x_i}{G}$ for $i\in [s]$. By assumption, $\eqt$ has $s$ 
blocks on $H$. Let $\blokt{y_1}{H}$, $\blokt{y_2}{H},\, \dotsc,$ 
$\blokt{y_s}{H}$ be the sequence of these blocks, where 
$y_1<y_2<\dotsb<y_s$ and $y_i=\max\blokt{y_i}{H}$ for every $i\in[s]$. 
By Lemma~\ref{lem-child}, the nodes $G$ and $H$ have $s$ children in 
$\calT^{M}$ and $\calT^{C}$. Let $G_i$ be the graph obtained from $G$ by 
addition of the edge $\{x_i,k+1\}$, let $H_i$ be the graph obtained from 
$H$ by addition of the edge $\{y_i,k+1\}$, for $i\in[s]$. By 
Lemma~\ref{lem-child}, the graphs $\{G_i;\ i\in[s]\}$ (or $\{H_i;\ 
i\in[s]\}$) are exactly the children of $G$ (or $H$, respectively). We 
define $\phi(G_i)=H_i$. The $\eqs\!$-blocks of $G_i$ are exactly the 
sets $\bloks{x_1}{G}$, $\bloks{x_2}{G},\,\dotsc,$ $\bloks{x_{i-1}}{G}$ 
and $\left(\bigcup_{j\ge i} \bloks{x_j}{G}\right)\setminus \{x_i\}$, 
while the $\eqt\!$-blocks of $H_i$ are exactly the sets 
$\blokt{y_1}{H}$, $\blokt{y_2}{H},\,\dotsc,$ $\blokt{y_{i-1}}{H}$ and 
$\left(\bigcup_{j\ge i} \blokt{y_j}{H}\right)\setminus \{y_i\}$. This 
implies that the corresponding blocks of $G_i$ and $H_i$ have the same 
number and the same size, as required (note that if $i=s$ and 
$\bloks{x_s}{G}=\{x_s\}$, then the last block in the above list of 
$\eqs\!$-blocks of $G_i$ is empty; however, this happens 
if and only if the last entry in the list of $\eqt\!$-blocks of $H_i$ is 
empty as well, so it does not violate the required properties of 
$\phi$). 

This concludes the proof of the theorem.
\end{proof}

\begin{cor}\label{cor-calc}
$M_{132}\eqp\calC$.
\end{cor}
\begin{proof}
Since $g(m,w,M_{132})$ is equal to the number of leaves of the tree 
$\calT^{M}$, and $g(m,w,\calC)$ is equal to the number of leaves of the 
tree $\calT^{C}$, this is a direct consequence of Theorem~\ref{thm-m132}.
\end{proof}

\begin{cor}\label{cor-sym}
$M_{132}\eqp M_{213}$. 
\end{cor}
\begin{proof}
Let $\overline{w}$ denote the Dyck word defined by the relation 
$\overline{w}_i=0$ if and only if $w_{n-i+1}=1$. By inverting the linear 
order of the vertices of a matching $M$ with base 
$w$, we obtain a matching $\overline{M}$ with base 
$\overline{w}$. Since every matching $C_k\in\calC$ satisfies 
$\overline{C_k}=C_k$, we know that a matching $M$ avoids $\calC$ if and 
only if $\overline{M}$ avoids $\calC$, and hence 
$g(m,w,\calC)=g(m,\overline{w}, \calC)$. Note that 
$\overline{M_{213}}=M_{132}$. This gives
\[
g(m,w,M_{132})=g(m,w,\calC)=g(m,\overline{w},\calC)=g(m,\overline{w},M_{132})=g(m,w,M_{213}),
\]
as claimed.
\end{proof}

\begin{cor}\label{cor-123-132} $M_{123}\gp M_{132}$. \end{cor}
\begin{proof}
Notice that $M_{123}=C_3\in \calC$, and all the other graphs in $\calC$ 
avoid $M_{123}$. This implies that 
$\calG(m,w,\calC)\subseteq\calG(m,w,M_{123})$, and for every $m\ge 4$ 
there is a $w\in\calD'_m$ for which this is a proper inclusion, because 
$C_m$ clearly belongs to $\calG(m,M_{123})\setminus\calG(m,\calC)$. The 
claim follows, as a consequence of Corollary~\ref{cor-calc}.
\end{proof}

\section{The matching $M_{231}$ and non-crossing pairs of Dyck paths}\label{sec-231}  

In this section, we prove that $M_{231}\eqp M_{321}$. This 
result is a special case of a more general theorem by Backelin et 
al.~\cite{bwx}. The new proof we present here provides a simple bijection 
between $M_{231}$-avoiding matchings of a fixed base $w$ and Dyck paths not 
exceeding the path representing the Dyck word $w$. 

We first introduce some notation: recall that 
$\calD(m)$ denotes the set of all Dyck paths of length $2m$. For two Dyck 
paths $P_1$ and $P_2$ of length $2m$, we say that $(P_1,P_2)$ is a 
\emph{non-crossing pair} if $P_2$ never reaches above $P_1$. Let $\calD^2_m$ 
denote the set of all the non-crossing pairs of Dyck paths of length $2m$ 
and, for a Dyck word $w$ of length $2m$, let $\calD^2_m(w)$ be the set of all 
the pairs $(P_1,P_2)\in\calD^2_m$ whose first component $P_1$ is the path 
represented by the Dyck word 
$w$.                                                      

Recently, Chen, Deng and Du \cite{chdd} have proved that $M_{123}\eqp M_{321}$ by 
a bijective construction involving Dyck paths. Their proof in fact shows 
that the cardinality of the set $D_m^2(w)$ is equal to the number of matchings 
with base $w$ avoiding $M_{123}$, and at the same time equal to the number of 
matchings with base $w$ avoiding $M_{321}$. In our notation, this corresponds 
to the following claim:
\[
\forall m\in\Nset\ \forall w\in\calD'_m\  g(m,w,M_{123})=|D^2_m(w)|=g(m,w,M_{321}).
\]
                                                                              
In this section, we extend these equalities to the matching $M_{231}$ by 
proving $|D^2_m(w)|=g(m,w,M_{231})$. This shows that $M_{231}\eqp M_{321}\eqp 
M_{123}$.
                                     
We remark that the number of non-crossing pairs of Dyck paths of length $2m$ 
(and hence the number of $M$-avoiding matchings of size $m$, where $M$ is any 
of $M_{123},M_{321}$ or $M_{231}$) is equal to $c_{m+2}c_m-c_{m+1}^2$, where 
$c_m$ is the $m$-th Catalan number (see \cite{gobe}). The sequence 
$(c_{m+2}c_m-c_{m+1}^2;\ m\in\Nset)$ is listed as the entry A005700 in the 
On-Line Encyclopedia of Integer Sequences \cite{oeis}. It is noteworthy, that 
N. Bonichon \cite{boni} has shown a completely different combinatorial 
interpretation of this sequence, in terms of realizers of plane 
triangulations.

Let us fix $m\in\Nset$ and $w\in\calD'_m$. Let $M$ be a matching with base 
$w$. Let $1=x_1<\dotsb<x_m$ denote the sequence of all the l-vertices with 
respect to $w$, and $y_1<\dotsb<y_m=2m$ be the sequence of all the 
r-vertices with respect to $w$. Let $y_k$ be the neighbour of 
$x_1$ in $M$. An edge $\{x_i,y_j\}$ of $M$ is called \emph{short} if 
$y_j<y_k$, and it is called \emph{long} if $y_j>y_k$. Let $E_\text{S}(M)$ and $E_\text{L}(M)$ 
denote the set of the short edges and long edges, respectively, so that we 
have $E(M)=E_\text{S}(M)\cup E_\text{L}(M)\cup \bigl\{\{1,y_k\}\bigr\}$. An l-vertex $x_i$ 
is called \emph{short} (or \emph{long}) if it is incident with a short edge 
(or a long edge, respectively).

\begin{lem}
Let $M$ be a matching with base $w$. $M$ avoids $M_{231}$ if 
and only if $M$ satisfies the following three conditions:
\begin{itemize}
\item
The subgraph of $M$ induced by the short edges avoids $M_{231}$.
\item
The subgraph of $M$ induced by the long edges avoids $M_{231}$.
\item
Every short l-vertex precedes all the long l-vertices.
\end{itemize}
\end{lem}
\begin{proof}
The first two conditions are clearly necessary. The third condition is 
necessary as well, for if $M$ contained an edge $\{x_s,y_s\}\in E_S$ and 
an edge $\{x_l,y_l\}\in E_L$ with $x_s>x_l$, then the six vertices 
$1<x_l<x_s<y_s<y_k<y_l$ would induce a copy of $M_{231}$. 

To see that the three conditions are sufficient, assume for 
contradiction that a graph $M$ satisfies these conditions but contains 
the forbidden configuration induced by some vertices 
$x_a<x_b<x_c<y_d<y_e<y_f$. We first note that $y_f>y_k$: indeed, it is 
impossible to have $y_f=y_k$, because $y_f$ is not connected to the 
leftmost vertex, and the inequality $y_f<y_k$ would imply that all the 
three edges of the forbidden configuration are short, which is ruled out 
by the first condition of the lemma. Thus, the edge $\{x_b,y_f\}$ is 
long, and hence $\{x_c,y_d\}$ is long as well, by the third condition. 
This implies that $y_d>y_k$, hence $y_e>y_k$ as well, and all the three 
edges of the configuration are long, contradicting the second condition 
of the lemma.
\end{proof}
            
\begin{figure}[ht]
\hfil
\includegraphics{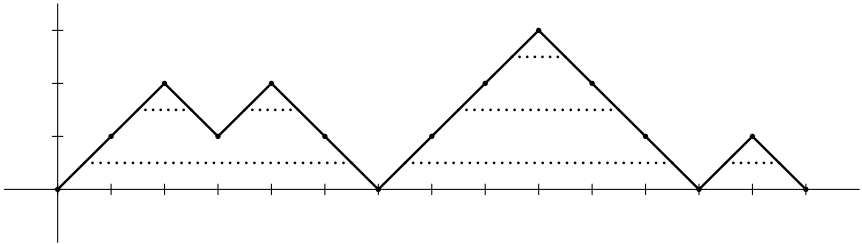}
\caption[A Dyck path with tunnels]{An example of a Dyck path. The dotted 
segments represent the tunnels.}\label{fig-tunnel}
\end{figure}

To construct the required bijection between $\calG(m,w,M_{231})$ and 
$\calD^2_m(w)$, we will use the intuitive notion of a ``tunnel'' in a Dyck 
path, which has been employed in bijective constructions involving 
permutations in, e.g., \cite{eliza1} or \cite{eliza2}. Let $P$ be a Dyck 
path. A \emph{tunnel} in $P$ is a horizontal segment $t$ whose left endpoint 
is the center of an up-step of $P$, its right endpoint is the center of a 
down-step of $P$, and no other point of $t$ belongs to $P$ (see 
Fig.~\ref{fig-tunnel}). A path $P\in\calD_m$ has exactly $m$ tunnels. An 
up-step $u$ and a down-step $d$ of $P$ are called \emph{partners} if $P$ has 
a tunnel connecting $u$ and $d$. Let $u_1(P),\dotsc,u_m(P)$ denote the 
up-steps of $P$ and $d_1(P),\dotsc,d_m(P)$ denote the down-steps of $P$, in 
the left-to-right order in which they appear on $P$.

\begin{figure}[ht]
\hfil\includegraphics[scale=.8]{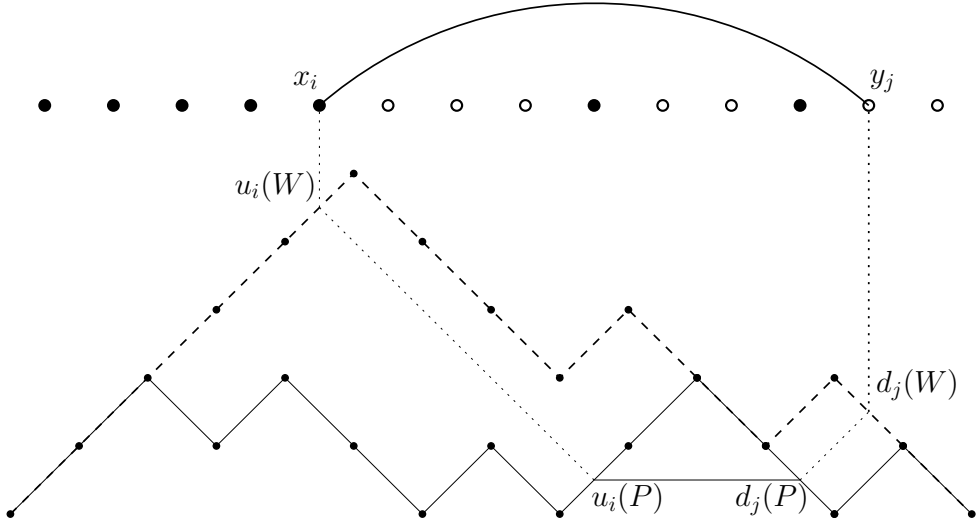} \caption[The correspondence between 
tunnels and edges]{The correspondence between a pair of Dyck paths $(W,P)$ 
and a matching $M(W,P)$. A tunnel between $u_i(P)$ and $d_j(P)$ 
corresponds to an edge $\{x_i,y_j\}$. The filled dots above the pair of paths 
represent the l-vertices of the matching, the empty dots represent the r-vertices.}\label{fig-tuhr}
\end{figure}

Let $W\in\calD_m$ be the Dyck path represented by the Dyck word $w$, let 
$(W,P)\in\calD^2_m(w)$ be a non-crossing pair of Dyck paths. Let $M(W,P)$ be 
the unique matching with base $w$ satisfying 
the condition that $\{x_i,y_j\}$ is an 
edge of $M$ if and only if $u_i(P)$ is the partner of $d_j(P)$. 
To see that this definition is valid, we need to check that 
if $u_i(P)$ is partnered to $d_j(P)$ 
in a path $P$ not exceeding $W$, then $x_i<y_j$ in the matchings with base 
$w$. This is indeed the case, because the horizontal coordinate of $u_i(W)$ 
(which  determines the position of $x_i$ in the matching) does not exceed the 
horizontal coordinate of $u_i(P)$, while the horizontal coordinate of $d_j(P)$ 
does not exceed the horizontal coordinate of $d_j(W)$ (note that a half-line 
starting in the center of $u_i(P)$ directed north-west intersects $W$ in the 
center of $u_i(W)$; similarly, a half-line starting in the center of $d_j(P)$ 
directed north-east hits the center of $d_j(W)$). See Fig.~\ref{fig-tuhr}.

\begin{lem}\label{lem-gpm231}
If $(W,P)\in\calD^2_m(w)$, then $M(W,P)$ avoids $M_{231}$.
\end{lem}
\begin{proof}
Choose an arbitrary $(W,P)\in\calD^2_m(w)$ and assume, for 
contradiction, that there are six vertices $x_a<x_b<x_c<y_d<y_e<y_f$ in 
$M(W,P)$ which induce the forbidden configuration. Let $t_{cd}, t_{ae}$
and $t_{bf}$ be the tunnels corresponding to the three edges $x_cy_d$, 
$x_ay_e$ and $x_by_f$, respectively. Note that the projection of $t_{cd}$
onto some horizontal line $h$ is a subset of the projections of 
$t_{ae}$ and $t_{bf}$ onto $h$. 
Thus, the three tunnels lie on different horizontal 
lines and there is a vertical line intersecting all of them.  

Since $a<b$, the tunnel $t_{ae}$ must lie below $t_{bf}$, otherwise the subpath
of $P$ between $u_a(P)$ and $u_b(P)$ would intersect $t_{ae}$.
On the other hand, $e<f$ implies that $t_{ae}$ lies above $t_{bf}$, 
a contradiction. 
\end{proof}

The aim of the next lemma is to show that the mapping $P\mapsto M(W,P)$ can 
be inverted.

\begin{lem}\label{lem-indm231}
For every $M\in\calG(m,w,M_{231})$ there is a unique Dyck path $P$ such 
that $(W,P)\in\calD_m^2(w)$ and $M=M(W,P)$.
\end{lem}
\begin{proof}
We proceed by induction on $m$. The case $m=1$ is clear, so let us assume 
that $m>1$ and that the lemma holds for every $m'<m$ and every 
$w'\in\calD'_{m'}$. Let us choose an arbitrary $w\in\calD'_m$, and an 
arbitrary $M\in\calG(m,w,M_{231})$, and define $k$ such that $\{1,y_k\}$ is 
an edge of $M$. Let $M_\text{S}$ be the matching from $\calG(k-1)$ that is 
isomorphic to the subgraph of $M$ induced by the short edges, let 
$M_\text{L}\in\calG(m-k)$ be isomorphic to the subgraph induced by the long 
ones, let $w_\text{S}$ and $w_\text{L}$ be the respective bases of 
$M_\text{S}$ and $M_\text{L}$, and let $W_\text{S}$ and $W_\text{L}$ be the 
Dyck paths corresponding to $w_\text{S}$ and $w_\text{L}$. By induction, we 
know that $M_\text{S}=M(W_\text{S},P_\text{S})$ and 
$M_\text{L}=M(W_\text{L},P_\text{L})$ for some Dyck paths $P_\text{S}$ and 
$P_\text{L}$, where $P_\text{S}$ does not exceed $W_\text{S}$, and 
$P_\text{L}$ does not exceed $W_\text{L}$. Let $w_\text{X}$ be the Dyck word 
$0w_\text{S}1w_\text{L}$, and let $W_\text{X}$ be the corresponding Dyck path. Note 
that $W_\text{X}$ does not exceed $W$: assume that $W$ has $t$ up-steps occurring 
before the $k$-th down-step; then $W_\text{X}$ is obtained from $W$ by omitting the 
$t-k$ up-steps $u_{k+1}(W), u_{k+2}(W),\dotsc,u_t(W)$, and inserting $t-k$ 
new up-steps directly after the $k$-th down-step. 

Let $P$ be the Dyck path obtained by concatenating the 
following pieces:
\begin{itemize}
\item
An up-step from $(0,0)$ to $(1,1)$
\item
A shifted copy of $P_\text{S}$ from $(1,1)$ to $(2k-1,1)$
\item
A down-step from $(2k-1,1)$ to $(2k,0)$
\item
A shifted copy of $P_\text{L}$ from $(2k,0)$ to $(2m,0)$
\end{itemize}
Since $P$ clearly does not exceed $W_\text{X}$, it does not exceed $W$ either. 
Let us check that $M=M(W,P)$:
\begin{itemize}
\item The base of $M$ is equal to the base of $M(W,P)$. Thus, to see that 
$M$ is equal to $M(W,P)$, it suffices to check that $M_\text{S}$ and $M_\text{L}$ are isomorphic to the 
matchings induced by the short edges of $M(W,P)$ and the long edges of 
$M(W,P)$, respectively.
\item The up-step $u_1(P)$ is clearly partnered to the down-step $d_k(P)$ 
(which connects $(2k-1,1)$ to $(2k,0)$). Thus, $M(W,P)$ contains the edge 
$\{x_1,y_k\}$. It follows that $M(W,P)$ has $k-1$ short edges, incident to 
the l-vertices $x_2,\dotsc,x_k$ and r-vertices $y_1,\dotsc,y_{k-1}$.
\item The $k-1$ up-steps $u_2(P),\dotsc,u_k(P)$ as well as the 
$k-1$ down\nobreakdash-steps $d_1(P),\dotsc,d_{k-1}(P)$ all belong to the shifted copy 
of $P_\text{S}$. Since 
shifting does not affect the partnership relations, we see that the short 
edges of $M(W,P)$ form a matching isomorphic to $M_\text{S}=M(W_\text{S},P_\text{S})$.  
\item Similarly, the up-steps $u_{k+1}(P),$ $u_{k+2}(P),\dotsc,u_m(P)$ are partnered to 
the down-steps $d_{k+1}(P)$, $d_{k+2}(P),\dotsc,d_m(P)$ according to the tunnels of $P_\text{L}$. The 
corresponding long edges form a matching isomorphic to~$M_\text{L}$.
\end{itemize}                                                
It follows that $M=M(W,P)$.

We now show that $P$ is determined uniquely: assume that $M=M(W,Q)$ for some 
$Q\in\calD_m$. Since $\{1,y_k\}\in E(M)$, the path $Q$ must contain a down 
step from $(2k-1,1)$ to $(2k,0)$, and this down-step must be the first 
down-step of $Q$ to reach the line $y=0$. This shows that the subpath of $Q$ 
between $(1,1)$ and $(2k-1,1)$ is a shifted copy of some Dyck path 
$Q_S\in\calD_{k-1}$. The tunnels of this path must define a matching 
isomorphic to $M_\text{S}=(W_\text{S},P_\text{S})$. By induction, we know 
that $P_\text{S}$ is determined uniquely, hence $P_\text{S}=Q_\text{S}$. By the same 
argument, we see that the subpath of $Q$ from $(2k,0)$ to $(2m,0)$ is a 
shifted copy of $P_\text{L}$. This shows that $P=Q$, and $P$ is unique, as 
claimed.
\end{proof}

We are now ready to prove the following theorem:

\begin{thm}\label{thm-m231}
For each $m\in\Nset$ and for each $w\in\calD'_m$, $g(m,w,M_{231})$ is 
equal to $|\calD^2_m(w)|$.
\end{thm}
\begin{proof}
Putting together Lemma~\ref{lem-gpm231} and Lemma~\ref{lem-indm231}, we infer 
that the function that maps a pair $(W,P)\in\calD_m(w)$ to the matching $M(W,P)$ 
is a bijection between $\calD^2_m(w)$ and $\calG(m,w,M_{231})$. This gives the 
required result.
\end{proof} 

\begin{cor}
$M_{231}\eqp M_{321}$.
\end{cor}
\begin{proof}
This a direct consequence of Theorem~\ref{thm-m231} and the results of Chen, Deng and Du \cite{chdd}.
\end{proof}

\section{Conclusion and Open Problems}
We have introduced an equivalence relation $\eqp$ on the set of 
permutational matchings, and we have determined how this equivalence partitions 
the set of permutational matchings of order 3. However, many natural 
questions remain unanswered. For instance,
it would be nice to have an estimate on the number of blocks of $\eqp$ belonging to the 
set $\calG(m)$. Also, is it possible to characterize the minimal and the maximal 
elements of $(\calG(m),\leqp)$?

\section*{Acknowledgements}
I am grateful to Martin Klazar for his comments and suggestions. I 
would also like to thank Anna de Mier and Zvezdelina Stankova for pointing out 
useful references to related research.

\end{document}